\documentclass[a4paper,12pt]{article}
  \usepackage[latin1]{inputenc}
  \usepackage[T1]{fontenc}
  \usepackage{amsmath,amssymb,epsfig,graphicx}
  \usepackage[all]{xy}

\def\R{\mathbb{R}}
\def\C{\mathbb C}

\def\P{\mathbb P}
\def\K{\mathbf{K}}
\title{Cyclicity of period annuli and principalization of Bautin ideals}
 \author{Lubomir Gavrilov \\
 \normalsize \it Institut de Math\'{e}matiques de Toulouse, UMR 5219\\
 \normalsize \it Universit\'{e} Paul-Sabatier (Toulouse III)\\  \normalsize \it 31062 Toulouse, Cedex 9\\ \normalsize \it   France  }
\begin{document}
\maketitle
\newtheorem{definition}{Definition}
\newtheorem{remark}{Remark}
\newtheorem{theorem}{Theorem} 
\newtheorem{lemma}{Lemma}
\newtheorem{proposition}{Proposition}
\newtheorem{corollary}{Corollary}
\vspace{5mm}
\begin{abstract}
Let $\Pi$ be an open period annulus of a plane analytic vector field $X_0$. We prove that the maximal number of
limit cycles which bifurcate from $\Pi$ under a given multi-parameter analytic deformation $X_\lambda$ of $X_0$
is the same as in an appropriate one-parameter analytic deformation $X_{\lambda(\varepsilon)}$, provided that
this cyclicity is finite. Along the same lines we give also a bound of the cyclicity of homoclinic saddle loops.

\end{abstract}
\section{Statement of the result}
Let $X_\lambda$, $\lambda\in (\R^n,0$) be an analytic family of plane vector fields, and let
 $\Pi$ be an open period annulus of $X_0$ (an open domain which
is a union of periodic orbits of $X_0$).  The \emph{cyclicity} $Cycl(\Pi,X_\lambda)$ of  $\Pi$ with respect to
the deformation $X_\lambda$ is the maximal number of limit cycles of $X_\lambda$ which tend to $\Pi$ as
$\lambda$ tends to zero, see Definition \ref{cyclicity} bellow. The reader should not confuse the cyclicity
$Cycl(\Pi,X_\lambda)$ of the open period annulus $\Pi$ with the cyclicity $Cycl(\bar{\Pi},X_\lambda)$ of the
closed period annulus $\bar{\Pi}$.

Our first result is the following
\begin{theorem}
\label{main}
If the cyclicity $Cycl(\Pi,X_\lambda)$ of the open period annulus $\Pi$ is finite, then there
exists a germ of analytic curve $\varepsilon \mapsto \lambda(\varepsilon), \varepsilon \in (\mathbb{R},0)$,
$\lambda(0)=0$,
 such that
\begin{equation}
\label{main1}
Cycl(\Pi,X_\lambda) = Cycl(\Pi,X_{\lambda(\epsilon))}).
\end{equation}
\end{theorem}
In other words, the problem of finding the cyclicity of an open period annulus with resect to a multi-parameter
deformation, can be always reduced to the "simpler" problem of finding cyclicity with respect to a one-parameter
deformation. Indeed, in this case the displacement map can be expanded into power series
\begin{equation}\label{ppm}
d(u,\epsilon)= \epsilon^k (M_\xi(u) + \epsilon R(u,\epsilon))
\end{equation}
where $M_\xi$, the so called \emph{higher order Poincar\'{e}-Pontryagin (or Melnikov) function}, depends on the germ
of analytic curve $\xi : \varepsilon \mapsto \lambda(\varepsilon)$. If $\Delta=(0,1)$ is an interval
parameterizing a cross section to the annulus $\Pi$, then $M_\xi$ is analytic on $(0,1)$ and its number of zeros
$Z(M_\xi)$ counted with multiplicity is an upper bound of $Cycl(\Pi,X_{\lambda(\epsilon))})$. Therefore we get
the inequality
\begin{equation}\label{main2}
Cycl(\Pi,X_\lambda) \leq \sup_{\xi} Z(M_\xi)
\end{equation}
where the $\sup$ is taken along all  germs of analytic curves $$\xi : \varepsilon \mapsto \lambda(\varepsilon),
\xi(0)=0 .
$$

Let $X_\lambda$ be an arbitrary analytic deformation of a vector field $X_0$. It has been conjectured by
Roussarie \cite[p.23]{rous98} that the cyclicity $Cycl(\Gamma,X_\lambda)$ of every compact invariant set
$\Gamma$ of $X_0$ is finite. A particular case of this Conjecture is therefore that $Cycl(\Pi,X_\lambda)<
\infty$, that is to say the claim of Theorem \ref{main} holds without the finite cyclicity assumption on the
open period annulus. It follows from the proof of Theorem \ref{main} that if $Cycl(\Pi,X_\lambda) = \infty$,
then there exists a germ of analytic curve $\varepsilon \mapsto \lambda(\varepsilon), \varepsilon \in
(\mathbb{R},0)$, $\lambda(0)=0$, such that the corresponding higher order Poincar\'{e}-Pontryagin function
$M_\xi(u)$ defined by (\ref{ppm}) has an infinite number of zeros in the interval $\Delta$. In the particular
case when $X_0$ has an analytic first integral in a neighborhood of the closed period annulus $\bar{\Pi}$, the
analytic properties of $M_\xi(u)$ are studied in \cite{gav1,gav}. Using this, it  might be shown that $M_\xi(u)$
has a finite number of zeros on the open interval $\Delta$.

\begin{figure}
\begin{center}
\includegraphics[scale=0.5,angle=90]{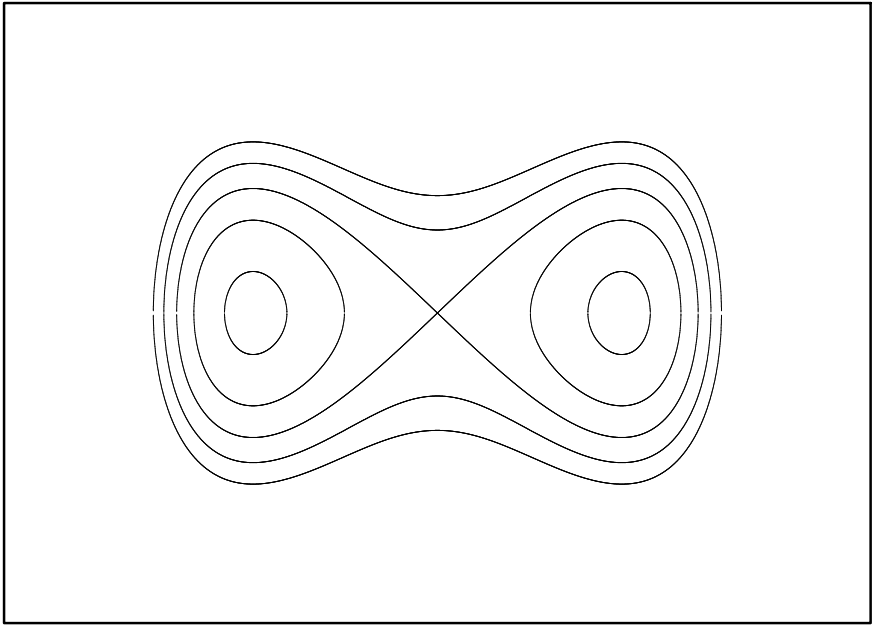}
\end{center}
\caption{Period annuli } \label{fig1}
\end{figure}

Suppose that the open period annulus of $X_0$  contains in its closure a non-degenerate center. Denote the union
of $\Pi$ with such a center by $\tilde{\Pi}$. If we blow up the center, it becomes a periodic orbit of a new
vector filed to which Theorem \ref{main} applies with minor modifications. It follows that (\ref{main1})
 holds true with $\Pi$ replaced by $\tilde{\Pi}$. The question whether we can replace the open period annulus
 $\Pi$ by its closure
$\overline{\Pi}$ is much more delicate. Namely, suppose that the closed period annulus $\overline{{\Pi}}$ is a
union of $\Pi$, a non-degenerate center, and a homoclinic saddle connection (as for instance the two bounded
annuli on fig.\ref{fig1}). By a homoclinic saddle connection (or separatrix loop) we mean a  union of
\emph{hyperbolic} saddle point with its stable and unstable separatrices which coincide. The union of these
separatrices is a homoclinic orbit of the vector field, and together with the saddle point they form a
separatrix loop or a homoclinic saddle connection.

 Suppose  that $\overline{\Delta}=[0,1]$ parameterizes
a cross section to $\overline{{\Pi}}$. Consider as above a germ of analytic curve $\xi$ and the corresponding
higher Poincar\'{e}-Pontryagin function $M_\xi$. It is continuous on $[0,1]$, analytic on $[0,1)$ and has an
asymptotic Dulac series at $u=1$
\begin{equation}\label{dulac}
M_\xi(u)= \sum_{i=0}^\infty a_{2i} (1-u)^i + a_{2i+1} (1-u)^{i+1} \log(1-u) \; .
\end{equation}
We define a {\it  generalized multiplicity  } of a zero of $M_\xi$ at $u=1$ to be equal to $j$ if $a_0=a_1=\dots
=a_{j-1}=0$, $a_j\neq 0$. At $u=0$ we define the multiplicity of the zero of
$$
M_\xi(u)= \sum_{i=1}^\infty a_{i} u^i
$$
to be equal to $j$, where $a_1=\dots =a_{j}=0$, $a_{j+1}\neq 0$.
 Define finally $ Z(M_\xi)$ to be the number of the zeros (counted with multiplicity)
of $M_\xi$ on $[0,1] $. It is classically known \cite{rous86} that
$$
Cycl(\overline{\Pi},X_{\lambda(\epsilon))}) \leq  Z(M_\xi) .
$$
We shall prove
\begin{theorem}
\label{main3} Under the above conditions
\begin{equation}\label{main4}
Cycl(\overline{\Pi},X_\lambda) \leq \sup_{\xi} Z(M_\xi)
\end{equation}
where the upper bound is taken along all  germs of analytic curves $$\xi : \varepsilon \mapsto
\lambda(\varepsilon), \xi(0)=0 .$$
\end{theorem}
\noindent \textbf{Remarque.} \emph{The hyperbolicity of the saddle point is essential for the proof of Theorem
\ref{main3}. It follows from \cite[Proposition A 2.1, p. 111]{rous89} that the vector field $X_0$ possesses an
analytic first integral in a neighborhood of $\overline{{\Pi}}$. Therefore $X_0$ is a Hamiltonian vector field
with respect to a suitable area form. It is also known that the cyclicity of the closed period annulus in this
case is finite \cite{rous98}.}\\
\noindent \textbf{Example.} \emph{Theorem \ref{main} applies to the three open annuli on fig.\ref{fig1}, while
Theorem \ref{main3} applies only to the two closed and bounded period annuli on fig.\ref{fig1}.}

The proofs of Theorems \ref{main} and \ref{main3} are inspired by the Roussarie's paper \cite{rous01}, in which
the inequality (\ref{main4}) is shown to be true for a single regular periodic orbit. The equality (\ref{main1})
for a single regular periodic orbit is announced by Cauberg \cite{ca01} (see also \cite{ca06} for related
results). The finiteness of the cyclicity $Cycl(\gamma,X_\lambda)$ of a regular periodic orbit $\gamma$ was
previously proved by Fran\c{c}oise and Pugh \cite{frpu}. The main technical tools in the proof of Roussarie's
theorem \cite{rous01} are the Hironaka's desingularization theorem applied to the Bautin ideal, followed by a
derivation-division algorithm. This second argument applies only locally. To prove Theorem \ref{main} we also
use the Hironaka's theorem, but replace the derivation-division algorithm by a variant of  the Weierstrass
preparation theorem. This already gives an upper bound of the cyclicity in terms of zeros of Poincar\'{e}-Pontryagin
functions. To get the exact result (\ref{main1}) we use the curve selection Lemma, as suggested by Roussarie
\cite{rous01}. Theorem \ref{main3} has a similar proof but is based on \cite[Theorem C]{rous86}. We note that we
obtain a non-necessarily exact upper bound for the cyclicity of the closed period annulus. The reason is that
the bifurcation diagram of limit cycles near the separatrix loop is not analytic, and therefore we can not apply
the curve selection lemma.

We mention finally that Theorem \ref{main} end \ref{main3} allow   an obvious complex version (with the same
proof). For this reason most of the results in the next section are stated in a complex domain as well. The
question about finding an explicit  upper bound of the cyclicity of more complicated separatrix connections is
almost completely opened. A recent progress in this direction is obtained in \cite{dumo}.

The paper is organized as follows. In the next section we formulate several classical results which will be used
in the proof. The latter is given in section \ref{proof}. In the last section \ref{remarks} we  discuss some
open questions.

\section{Digression}
In this section we formulate, for convenience of the reader, several facts of general interest, which are
necessary for the proof of Theorem \ref{main}. The base field is $\mathbf{K}=\mathbb{R}$ or $\mathbb{C}$. The
corresponding projective space $\mathbb{P}_\mathbf{K}$ is denoted simply $\mathbb{P}$.

\subsection{Principalization of ideals}

 Let $\varphi_0,\varphi_1,\dots,\varphi_p$ be non-zero analytic functions on a smooth
complex or real analytic variety $X$.
The indeterminacy points
of the rational map
$$
\varphi: X \dashrightarrow \P^p
$$
can be eliminated as follows \cite{hiro,bier}

\begin{theorem}[Hironaka desingularization]
\label{hiroth}
 There exists a smooth analytic variety $\tilde{X}$ and a proper analytic map $\pi:
\tilde{X}\rightarrow X$ such that the induced map $\tilde{\varphi}=\varphi\circ \pi$ is analytic.
$$
 \xymatrix{
    \tilde{X}  \ar[d]_\pi \ar[rd]^{\tilde{\varphi}}&  \\
    X \ar@{-->}[r]^\varphi & \P^n }
  $$
  \end{theorem}
 Let $\mathcal{O}_X$ be the sheaf of analytic functions on $X$ and consider the ideal
sheaf $I \subset \mathcal{O}_X$ generated by  $\varphi_0,\varphi_1,\dots,\varphi_p$. The inverse image ideal
sheaf of $I$ under the map $\pi: \tilde{X}\rightarrow X$ will be denoted $\pi^*I$. This is the ideal sheaf
generated by the pull-backs of local sections of $I$. We note that $\pi^*I$  may differ from the usual
sheaf-theoretic pull-back, also commonly denoted by $\pi^*I$. A simple consequence of Theorem \ref{hiroth} is
the following

\begin{corollary}
\label{hirocor} The inverse image ideal sheaf $\pi^*I$ is principal.
 \end{corollary}
This is called the principalization of $I$. Indeed, as the induced map $\tilde{\varphi}$ is analytic, then for
every $\tilde{\lambda}\in \tilde{X}$ there exists  $j$, such that the functions $\tilde{\varphi}_i/
\tilde{\varphi}_j$, $i=1,2,\dots,p$, are analytic in a neighborhood of $\tilde{\lambda}$. Therefore there is a
neighborhood $\tilde{U}$ of $\tilde{\lambda}$ such that $\tilde{\varphi}_j|_{\tilde{U}}$ divides
$\tilde{\varphi}_i|_{\tilde{U}}$ in the ring of sections $\mathcal{O}_{\tilde{U}}$ of the sheaf
$\mathcal{O}_{\tilde{X}}$, that is to say $I_{\tilde{U}}$ is generated by $\tilde{\varphi}_j|_{\tilde{U}}$.
\subsection{The Weierstrass preparation Theorem}
\begin{definition}
Let $[a,b]\subset \R$. A Weierstrass polynomial in a  neighborhood of $[a,b]\times\{0\}\subset \K\times \K^n$ is
an analytic function of the form
$$
P(u,\lambda)= u^d+a_1(\lambda) u^{d-1}a_2(\lambda)+  \dots + a_d(\lambda)
$$
such that $P(u,0)$ has exactly $d$ zeros in $[a,b]$ (counted with multiplicity).
\end{definition}
In the case $a=b=0$, $\K=\C$, the above definition coincides with the usual definition of a Weierstrass
polynomial in a neighborhood of the origin in $\C^{n+1}$ \cite{gunn}.

\begin{theorem}[Weierstrass preparation theorem]
\label{weie}
 Let $f(u,\lambda)$ be an analytic function in a  neighborhood of $[a,b]\times\{0\}\subset
\K\times \K^n$ such that $f(u,0)$ is not identically zero. Then $f$ has an unique representation $ f= P.h$ where
$P=P(u,\lambda)$ is a Weierstrass polynomial in a neighborhood of $[a,b]\times\{0\}\subset \K\times \K^n$, and
$h=h(u,\lambda)$ is an analytic function, such that $h(u,0)\neq 0, \forall u \in [a,b]$.
\end{theorem}
The proof of the above theorem is  the same as in the usual case $a=b$, $\K=\C$ \cite{gunn}.

We are also interested in the behavior of the zeros $u=u(\lambda)$ of the Weierstrass polynomial. For this
reason we consider the discriminant $\Delta(\lambda)$ of $P(u,\lambda)$ with respect to $u$. It is an analytic
function in a neighborhood of the origin $0\in \K^n$ which might be also identically zero. This may happen for
instance if for every fixed $\lambda$, such that $\|\lambda\|$ is sufficiently small, the polynomial
$P(u,\lambda)$ has a double zero $u(\lambda)$, which is then analytic in $\lambda$. The analytic function
$(u-u(\lambda))^2$ then divides the Weierstrass polynomial $P(u,\lambda)$, the result being also a Weierstrass
polynomial. These considerations generalize to the following
\begin{corollary}
\label{weie1} The Weierstrass polynomial $P$ from Theorem \ref{weie} has a representation
$P=P_1^{i_1}.P_2^{i_2}\dots P_k^{i_k}$, where for each $i$, $P_i$ is a Weierstrass polynomial in a neighborhood
of $[a,b]\times\{0\}\subset \C\times \C^n$, with non-vanishing identically discriminant
$\Delta_i(\lambda)\not\equiv 0$.
\end{corollary}
It follows that the \emph{bifurcation locus} of the zeros $\{ u:  P(u,\lambda)=0\}$ of a Weierstrass polynomial
$P(.,\lambda)$ (and hence of an analytic function $f(.,\lambda)$) in a neighborhood of a compact interval
$[a,b]$ is a germ of an analytic set in a neighborhood of the origin in $\C^n$. Namely, in the complex case
$\K=\C$, this bifurcation locus $\mathcal{B}^\mathbb{^C}$ is the union of discriminant loci of $P_i(.,\lambda)$
and resultant loci of pairs $P_i(.,\lambda),P_j(.,\lambda)$.

In the real case  $\K=\C$ the real bifurcation locus $\mathcal{B}\mathbb{^R}$ is contained in the real part $\Re
(\mathcal{B}\mathbb{^C})$ of the complex one and is also a real analytic set.
 The complement to the bifurcation locus $\mathcal{B}\mathbb{^R}$ of real zeros of a  real analytic function
 $f(.,\lambda)$ is therefore  a \emph{semi-analytic set}. The more general case, when $f(u,0)\equiv 0$ on $[a,b]$
 will be considered in section \ref{proof}. It will follow from the proof of Theorem \ref{main} that the
 complement to the bifurcation locus $\mathcal{B}\mathbb{^R}$ is then a \emph{sub-analytic set}, see
 also Caubergh \cite{ca01}.

\subsection{The curve selection Lemma}
\begin{lemma}
Let $U$ be an open neighborhood of the origin in $\R^n$ and let $$f_1,\dots,f_k, g_1,\dots,g_s$$ be real
analytic functions on $U$ such that the origin is in the closure of the semi-analytic set:
$$ Z := \{ x\in  U: f_1(x) = \dots = f_k(x) \mbox{  and  } g_i(x) > 0, i = 1,\dots, s\}
$$
Then there exists a real analytic curve
 $\gamma : [0, \delta) \rightarrow U $ with
$\gamma(0) = 0$ and $\gamma(t)\in Z$, $ \forall t \in (0,\delta)$.
\end{lemma}
We refer to Milnor's book \cite{miln} for a proof. He does it in the algebraic category, but his proof works in
general with minor (obvious) modifications.

Let $f=f(u,\lambda)$ be a function, real analytic  in a real neighborhood of $[a,b]\times\{0\}\subset \R\times
\R^n$ (we do not suppose that $f(.,0)\not\equiv 0$).
\begin{definition}
\label{cyclf} The cyclicity $Cycl(([a,b],f(.,0)), f(.,\lambda))$ of $ f(.,0)$ on the interval $[a,b]\subset\R$
is the smallest integer $N$ having the property: there exists $\varepsilon_0 > 0$ and a real neighborhood $V$ of
$[a,b]$, such that for every $\lambda\in\R^n$, such that $\|\lambda\| < \varepsilon_0$, the function
$f(u,\lambda)$  has no more than $N$ distinct zeros in $V$, counted without multiplicity.
\end{definition}
When there is no danger of confusion  we shall write  $Cycl([a,b], f(.,\lambda))$ instead of
$Cycl(([a,b],f(.,0)), f(.,\lambda))$.
 The number of the zeros of $f(.,0)$ on the interval $[a,b]\subset
\R$, counted with multiplicity, is an upper bound for $Cycl([a,b], f(.,\lambda))$, but not necessarily an exact
bound. The Weierstrass preparation theorem, Corollary \ref{weie1} and the Curve selection Lemma imply
\begin{theorem}
\label{baby} Let $f(u,\lambda)$ be a real analytic function in a  neighborhood of $[a,b]\times\{0\}\subset
\R\times \R^n$ which is non-identically zero on $[a,b]\times\{0\}$. There exists an analytic curve $[0; \delta)
\rightarrow \R^n :\varepsilon \mapsto \lambda(\varepsilon)$, $\lambda(0)=0$, such that
$$
Cycl([a,b], f(.,\lambda)) = Cycl([a,b], f(.,\lambda(\varepsilon))) .
$$
\end{theorem}

\subsection{Cyclicity of period annuli and the Bautin ideal} In this section the base field is $\K=\R$.
\begin{definition}
\label{cyclicity}
 Let $X_\lambda$ be a family of analytic real plane vector fields depending analytically on a
parameter $\lambda \in (\R^n,0)$, and let $ K \subset R^2$ be a compact invariant set of $X_{\lambda_0}$. We say
that the pair $(K, X_{\lambda_0}$) has cyclicity $N = Cycl((K,X_{\lambda_0}), X_\lambda)$ with respect to the
deformation $X_\lambda$, provided that $N$ is the smallest integer having the property: there exists
$\varepsilon_0 > 0$ and a neighborhood $V_K$ of $K$, such that for every $\lambda$, such that $\|\lambda-
\lambda_0\| < \varepsilon_0$, the vector field $X_\lambda$  has no more than N limit cycles contained in $V_K$.
If $\tilde{K}$ is an invariant set of $X_{\lambda_0}$ (possibly non-compact), then the cyclicity of the pair
$(\tilde{K} , X_{\lambda_0})$ with respect to the deformation $X_\lambda$ is
$$
Cycl((\tilde{K} , X_{\lambda_0}),X_\lambda) = sup\{Cycl((K,X_{\lambda_0}), X_\lambda) : K \subset  \tilde{K} , K
\mbox{   is a compact }\}.
$$
\end{definition}
The cyclicity $Cycl((\tilde{K} , X_{\lambda_0}),X_\lambda)$ is therefore the maximal number of limit cycles
which tend to $\tilde{K}$ as $\lambda$ tends to $0$. To simplify the notation, and if there is no danger of
confusion, we shall write $ Cycl(K  ,X_\lambda) $ on the place of $ Cycl((K , X_{\lambda_0}),X_\lambda) . $
 Let $\Pi\subset \R^2$ be an open period annulus of a plane analytic vector field $X_0$. There
is a bi-analytic
 map identifying $\Pi$ to $S^1\times \Delta$ where $\Delta$ is a connected open interval.
 Therefore $X_0$ has an analytic first integral $u$ induced by the canonical projection
 $$ S^1\times \Delta \rightarrow \Delta: (\varphi,u)\mapsto u$$
 which parameterizes a
cross-section of the period annulus $\Pi$. Let $u \mapsto P(u,\lambda)$ be the first return map
 and $\delta(u,\lambda)= P(u,\lambda) - u$ the displacement function of $X_\lambda$. For every closed interval
 $[a,b]\subset \Delta$ there exists $\varepsilon_0 > 0$ such that the displacement function
 $\delta(u,\lambda)$ is well defined and analytic in
 $\{ (u,\lambda):    a-\varepsilon_0 < u < b+ \varepsilon_0, \| \lambda
 \|<\varepsilon_0\}$. For every fixed $\lambda$ there is a one-to-one correspondance between  zeros of
 $\delta(u,\lambda)$ and limit cycles of
 the vector field $X_\lambda$. This allows to define the cyclicity $Cycl(\Pi,X_\lambda)$ in terms of the
 cyclicity of the displacement function $\delta(u,\lambda)$ on the cross section $\Delta$ (Definition
 \ref{cyclf}):

\begin{equation}\label{c1}
Cycl(K,X_\lambda)= Cycl([a,b],\delta(.,\lambda))
\end{equation}
where $K= S^1
 \times [a,b]$ (we identified $\Pi$ and $S^1\times \Delta$) and

\begin{equation}\label{c2}
Cycl(\Pi,X_\lambda)= \sup_{[a,b]\subset \Delta}Cycl([a,b],\delta(.,\lambda)).
\end{equation}

 Let $u_0\in
\Delta$ and let us expand
$$
\delta(u,\lambda)= \sum_{i=0}^\infty a_i(\lambda) (u-u_0)^i.
$$
\begin{definition}[Bautin ideal \cite{rous89}, \cite{rous98}]
We define the Bautin ideal $\mathcal{I}$ of $X_\lambda$ to be the ideal generated by the germs $\tilde{a}_i$ of
$a_i$ in the local ring $\mathcal{O}_0(\R^n)$ of analytic germs of functions at $0\in \R^n$.
\end{definition}
This ideal is Noetherian and let $\tilde{\varphi}_1,\tilde{\varphi}_2,\dots,\tilde{\varphi}_p$ be a minimal
system of generators, where $p = \dim_\R \mathcal{I} / \mathcal{MI}$, and $\mathcal{M}$ is the maximal ideal of
the local ring $\mathcal{O}_0(\R^n)$. Let $\varphi_1,\varphi_2,\dots,\varphi_p$ be analytic functions
representing the generators  of the Bautin ideal in a neighborhood of the origin in $\R^n$.
\begin{proposition}[Roussarie, \cite{rous98}]
\label{bautin}
 The Bautin ideal
does not depend on the point $u_0\in \Delta$. For every $[a,b] \subset \Delta$ there is an open neighborhood $U$
of $[a,b]\times \{0\}$ in $\R\times \R^n$ and analytic functions $h_i(u,\lambda)$ in $U$, such that
\begin{equation}\label{b1}
\delta(u,\lambda)= \sum_{i=0}^p \varphi_i(\lambda) h_i(u,\lambda) .
\end{equation}
The real vector space generated by the functions $h_i(u,0), u\in [a,b]$ is of dimension $p$.
\end{proposition}
\section{Proof}
\label{proof}
To prove Theorem \ref{main} we may apply Theorem \ref{baby} to the displacement function $\delta(u,\lambda)$ ...
provided that $\delta(u,0)$ is not identically zero. This is certainly not the case. To overcome this difficulty
we principalize the Bautin ideal and divide the displacement map $\delta(u,\lambda)$ by a suitable analytic
function (which does not affect its cyclicity).

Suppose that the cyclicity $Cycl(\Pi,X_\lambda)$ is finite. There exists an invariant compact subset $K$ of
$X_0$ such that $Cycl(\Pi,X_\lambda)=Cycl(K,X_\lambda)$. If we identify $\Pi$ to $S^1\times \Delta$, then $K$ is
identified to $S^1 \times [a,b]$ where $[a,b]\subset \Delta$. It follows from (\ref{c1}), (\ref{c2}) and the
definition of cyclicity, that there is   an open interval $\sigma$, $[a,b] \subset \sigma\subset \Delta,$ and a
convergent sequence $(\lambda^k)_k$ in $\R^n$ which tends to the origin in $\R^n$, and such that for every $k$,
$\delta(.,\lambda^k)$ has exactly $ Cycl(\Pi,X_\lambda)$ zeros in the interval $\sigma$.

If $\varepsilon>0$ is sufficiently small, then the generators $\varphi_i$ of the Bautin ideal, defined by
(\ref{b1}), are analytic in the set $X=\{\lambda \in \C^n:\| \lambda \|<\varepsilon\}$. According to Theorem
\ref{hiroth} the rational map
$$
\varphi=(\varphi_0,\varphi_1,\dots,\varphi_n): X \dashrightarrow \P^n
$$
can be resolved. The projection $\pi: \tilde{X}\rightarrow X$ is a proper map which implies that there is a
convergent sequence $(\tilde{\lambda}^{k_i})_i$ in $\tilde{X}$, such that
$$
\pi(\tilde{\lambda}^{k_i})=\lambda^{k_i}, \lim_{i\rightarrow\infty} \tilde{\lambda}^{k_i} =
\tilde{\lambda}^{0}\in \pi^{-1}(0)
$$
and hence
\begin{eqnarray*}
Cycl(\Pi,X_\lambda) & = &  Cycl(S^1 \times [a,b],X_\lambda) \\
& = & Cycl((S^1 \times [a,b],X_{\pi(\tilde{\lambda}^0)}),X_{\pi(\tilde{\lambda})})\\
&= & Cycl( ([a,b],\delta(.,\pi(\tilde{\lambda}^0)),\delta(.,\pi(\tilde{\lambda}))).
\end{eqnarray*}
In other words, at $\tilde{\lambda}^0 \in \tilde{X}$, the cyclicity of the open period annulus is maximal.
Of course $X_{\pi(\tilde{\lambda}^0)} = X_0$, $\delta(.,\pi(\tilde{\lambda}^0)) = \delta(.,0)=0$.
Let $\tilde{\lambda}$ be a local variable on $\tilde{X}$ in a neighborhood of $\tilde{\lambda}_0$,
$\pi(\tilde{\lambda})=\lambda$. By Corollary \ref{hirocor} the inverse image of the Bautin ideal sheaf is
principal. Let $\varphi_0\circ\pi$ be a generator of the ideal of sections in a neighborhood of
$\tilde{\lambda}^0$. By Proposition \ref{bautin}, in a suitable neighborhood $U_{\tilde{\lambda}_0}$ of $
[a,b]\times \{\tilde{\lambda}_0\}$ in $\R\times \tilde{X}$ holds

\begin{equation}\label{principal}
\delta(u,\pi(\tilde{\lambda}))= \tilde{\varphi}_0(\tilde{\lambda}) \tilde{h}(u,\tilde{\lambda})
\end{equation}
where $$\tilde{\varphi}_0= \varphi_0\circ \pi, \tilde{h}(u,\tilde{\lambda})=h(u, \pi( \tilde{\lambda}))$$ and
$$h(u,
\pi( \tilde{\lambda}^0)= h(u,0) \not\equiv 0 .$$ We conclude that
$$
Cycl( ([a,b],\delta(.,\pi(\tilde{\lambda}^0)),\delta(.,\pi(\tilde{\lambda}))) = Cycl(
([a,b],h(.,\pi(\tilde{\lambda}^0)),h(.,\pi(\tilde{\lambda})))
$$
and by Theorem \ref{baby} there exists an analytic curve $\varepsilon \mapsto \tilde{\lambda}(\varepsilon)$,
$\tilde{\lambda}(0)= \tilde{\lambda}^0$, such that
$$
Cycl( ([a,b],h(.,\pi(\tilde{\lambda}^0)),h(.,\pi(\tilde{\lambda})))= Cycl(
([a,b],h(.,\pi(\tilde{\lambda}^0)),h(.,\pi(\tilde{\lambda}(\varepsilon)))).
$$
The curve $\varepsilon \mapsto \lambda(\varepsilon) =  \pi(\tilde{\lambda}(\varepsilon))$, $\lambda(0)=0$, is
analytic which shows finally that
$$
Cycl(\Pi,X_\lambda) = Cycl(\Pi,X_{\lambda(\varepsilon)}).
$$
Theorem \ref{main} is proved.$\Box$

The proof of Theorem \ref{main3} is similar: we resolve the map $\varphi$ and  principalize the Bautin ideal. Of
course we can not use the Weierstrass preparation theorem, neither the curve selection lemma. Let
$\tilde{\lambda}^0 \in \pi^{-1}(0)$ be a point in a neighborhood of which the closed period annulus has a
maximal cyclicity
$$
Cycl(\overline{\Pi},X_\lambda) = Cycl( (\overline{\Pi},X_{\pi(\tilde{\lambda}^0)}), X_{\pi(\tilde{\lambda})}).
$$
The displacement map takes the form  (\ref{principal}). Suppose that $\overline{\Delta}=[0,1]$ parameterizes a
cross section to $\overline{{\Pi}}$. Then $\tilde{h}(u,\tilde{\lambda_0})= h(u,0)\not\equiv 0$ is analytic on
$[0,1)$. First of all we have to prove that the Dulac expansion of $h(u,0)$ is not identically zero. Indeed, it
follows from \cite[Proposition A 2.1, p. 111]{rous89} that the vector field $X_\lambda$ possesses an analytic
first integral in a neighborhood of $\overline{{\Pi}}$. Therefore $X_\lambda$ is a Hamiltonian vector field with
respect to a suitable area form. Let
$$\widetilde{\xi} : \varepsilon \mapsto
\widetilde{\lambda}(\varepsilon), \widetilde{\lambda}(0)=\widetilde{\lambda}^0$$ be any analytic curve, not
contained in the zero locus of $\tilde{\varphi}_0(\tilde{\lambda})$. The displacement map of $X_\lambda$
restricted to the curve $\xi=\pi\circ \widetilde{\xi}$ is
\begin{equation}\label{eps}
\delta(u,\pi(\tilde{\lambda}(\varepsilon)))= c .\varepsilon ^k (h(u,0) + \varepsilon R(u,\varepsilon)), c\neq 0
\end{equation}
which shows that $h(u,0)=M_\xi(u)$ is the Poincar\'{e}-Pontryagin function associated to the curve $\xi$. It follows
from Roussarie's theorem  \cite[Theorem C]{rous89} that the Dulac expansion of $h(u,0)$ (\ref{dulac}) can not be
zero provided that $h(u,0)\not\equiv0$, and that the cyclicity of the loop $\gamma$ of $X_0$ with respect to the
deformation $X_{\pi(\tilde{\lambda}(\varepsilon))}$ is bounded by the generalized multiplicity of the zero $u=0$
of $h(u,0)$. Strictly speaking, the Roussarie's Theorem is proved for the case $k=1$ in (\ref{eps}), but the
proof in the case $k>1$ is exactly the same. The fact that $h(u,0)$ has a Dulac expansion (\ref{dulac}) follows
also from \cite{gav,gav1} where its monodromy is computed.

We conclude that the cyclicity of the closed period annulus is bounded by the total number of zeros $Z(M_\xi)$,
counted with generalized multiplicity. This completes the proof of Theorem \ref{main3}.$\Box$
\section{Concluding Remarks}
\label{remarks} Theorems \ref{main} and \ref{main3} seem to belong  to the mathematical folklore. The authors of
\cite{chow} for instance  used them to compute the cyclicity of open or closed period annuli of particular
quadratic systems with a center (but provided
 wrong references to \cite{ili98} and to the author's paper \cite{Gav4}, see \cite[Remark 2.1 and Lemma 2.1]{chow}).
Particular cases of  Theorems \ref{main} and \ref{main3} were previously used in \cite[p.223-224]{horil} and
 \cite[p.490-491]{Gav4}).

As in the Introduction, let $X_\lambda$ and $\Pi$ be an analytic family of analytic vector fields and an open
period annulus of the field $X_{\lambda_0}$ respectively. Let $\mathcal{M}_{X}=
\mathcal{M}_{X_\lambda}(X_{\lambda_0})$ be the set of all Poincar\'{e}-Pontryagin functions $M_\xi$ (\ref{ppm})
associated to germs of analytic curves $\xi : \varepsilon \mapsto \lambda(\varepsilon)$, $\lambda(0)=\lambda_0$.
The set $\mathcal{M}_{X}$ is not always a vector space, but spans a real vector space of finite dimension
(bounded by the number of generators of the Bautin ideal). Recall that all functions $M_\xi$ are defined on a
suitable open interval $\Delta$. We denote by $Z(M_\xi)$ the number of the zeros of $M_\xi$ on $\Delta$ (counted
with multiplicity). It follows from the proof of Theorem \ref{main} that
\begin{equation}\label{supm} \sup_{M_\xi \in
\mathcal{M}_{X}} Z(M_\xi) <\infty
\end{equation}if and only if for all $M_\xi \in \mathcal{M}_{X}$ holds
\begin{equation}\label{q1}
   Z(M_\xi)< \infty .
\end{equation}
As explained in the Introduction, the Roussarie's conjecture \cite[p.23]{rous98} (if it were true) would imply
that (\ref{supm}), (\ref{q1}) hold true.

Suppose now that $X_\lambda$, $\lambda \in \Lambda_n$, is the family of polynomial vector fields of degree at
most $n$ and denote
$$
Z(n,X_{\lambda_0})=\sup_{M_\xi \in  \mathcal{M}_{X}} Z(M_\xi) .
$$
According to Theorem \ref{main} the number $Z(n,X_{\lambda_0})$ is just the cyclicity of the open period annulus
$\Pi$
$$
Z(n,X_{\lambda_0})= Cycl(\Pi,X_\lambda) .
$$
If $X_{\lambda_0}$ is a generic Hamiltonian vector field with a center, the set $\mathcal{M}_{X}$ is a vector
space of Abelian integrals. The weakened  16th Hilbert problem, as stated by Arnold \cite{arn}, asks to compute
explicitly the number $Z(n,X_{\lambda_0})$. In this case the inequality (\ref{q1}) follows from the
Varchenko-Khovanskii theorem \cite{var,kho}.  Suppose now that $X_{\lambda_0}$ is a given plane vector field of
degree $n$ with a center (not necessarily Hamiltonian). The natural generalization of the Arnold's question is
then
\begin{center}
\emph{Find the numbers }$Z(n,X_{\lambda_0}), \lambda_0 \in \Lambda_n$.
\end{center}
To prove the finiteness of the number
$$
\sup_{\lambda_0\in \Lambda_n} Z(n,X_{\lambda_0})
$$
then would be a generalization of the Varchenko-Khovanskii theorem.

To answer the above questions, it is necessary to compute first the space of all Poincar\'{e}-Pontryagin functions
$$\mathcal{M}_{X}=
\mathcal{M}_{X_\lambda}(X_{\lambda_0}), \lambda \in \Lambda_n .$$ In the case $n=2$ this is a  result of Iliev
\cite[Theorem 2, Theorem 3]{ili98}. It is known that if $X_0$ is a quadratic Hamiltonian field, then
$Z(2,X_{0})=2$, except in the Hamiltonian triangle case, in which $Z(2,X_{0})=3$, see \cite{Gav4,chow}.

\emph{Acknowledgements}.  We acknowledge the critical comments of the referee, especially concerning the finite
cyclicity of the open period annuli.  We thank I.D. Iliev for the remarks which helped us to improve the text,
as well  D. Panazzolo for the explications concerning the proof of Theorem \ref{main3}.

\end{document}